\documentclass[12pt]{amsart}  
\usepackage{amsmath,amsthm,amsfonts,amssymb} 
\numberwithin{equation}{section} 
\newtheorem{thm}[equation]{Theorem}\newtheorem{cor}[equation]{Corollary} 
\newtheorem{lem}[equation]{Lemma} \newtheorem{prop}[equation]{Proposition} 
\newtheorem{example}[equation]{Example}
\theoremstyle{definition}  

\DeclareMathOperator{\Aut}{Aut}
\DeclareMathOperator{\codim}{codim}
\DeclareMathOperator{\diag}{diag}
\DeclareMathOperator{\End}{End}
\DeclareMathOperator{\Ext}{Ext}
\DeclareMathOperator{\Hom}{Hom}
\DeclareMathOperator{\id}{id}
\DeclareMathOperator{\Ima}{Im}
\DeclareMathOperator{\Ker}{Ker}
\DeclareMathOperator{\Span}{Span}

\newcommand{\C}{\mathbb{C}}
\newcommand{\DOT}{\setlength{\unitlength}{1pt}\begin{picture}(2.5,2)
                  (1,1)\put(2,3.5){\circle*{2}}\end{picture}}

\newcommand{\HH}{{\rm HH}}
\newcommand{\HHD}{{\rm HH}^{\DOT}}

\newcommand{\ot}{\otimes} 
\newcommand{\pe}{\perp}

\newcommand{\Wedge}{\textstyle\bigwedge}
\newcommand{\Z}{\mathbb {Z}}

\begin{document}

\title[SKEW DERIVATIONS AND DEFORMATIONS]
{SKEW DERIVATIONS AND DEFORMATIONS OF A FAMILY OF GROUP CROSSED PRODUCTS}

\author{Sarah Witherspoon}
\address{Department of Mathematics\\
 Texas A\&M University \\ College Station, Texas 77843, USA} \email{sjw@math.tamu.edu} 
\thanks{Partially supported by
NSF grants \#{}DMS--0422506 and \#{}DMS-0443476.
The author thanks the Mathematical Sciences Research
Institute for its hospitality during much of the writing of this paper.}

\date{October 31, 2005}

\begin{abstract}
We obtain deformations of a crossed product of a polynomial algebra
with a group, under some conditions, from universal deformation formulas.
We show that the resulting deformations are nontrivial by a comparison
with Hochschild cohomology.
The universal deformation formulas arise from actions of 
Hopf algebras generated by automorphisms and skew derivations, and
are universal in the sense that they apply to deform all algebras
with such Hopf algebra actions.
\end{abstract}

\maketitle

\section{Introduction}

Deformations of a polynomial algebra, such as the Weyl 
algebra or functions on quantum affine space, may be expressed 
by formulas
involving derivations of the polynomial algebra. 
These formulas are power series
in an indeterminate with coefficients in the universal enveloping
algebra of the Lie algebra of derivations. 
There are generalizations of such
deformations to other types of algebras, such as functions on a manifold
or orbifold, that are of current interest.

In this note we give a new generalization of the formulas themselves,
and apply it to crossed products of
polynomial algebras with groups of
linear automorphisms. These group crossed products 
are of interest in geometry due
to their relationship with corresponding orbifolds.
Particular deformations of such crossed products, called graded
Hecke algebras, were defined by Drinfel'd \cite{drinfeld86}.
These deformations have been studied by many authors, for example
for crossed products with real reflection groups, see \cite{lusztig89}, with
complex reflection groups, see \cite{ram-shepler},
and with symplectic reflection groups, see \cite{etingof-ginzburg}.
For these crossed product algebras, the universal enveloping
algebra of the Lie algebra of derivations does not
capture all the known deformations.
Instead we derive a deformation formula
from the action of a 
Hopf algebra under some hypotheses, recovering more of these
known deformations as well as some new ones.

We use the theory developed by
Giaquinto and Zhang of a universal deformation
formula based on a bialgebra $B$ \cite{giaquinto-zhang98}, extending
earlier such formulas based on universal enveloping algebras of Lie algebras. 
Such a formula is {\em universal} in the sense
that it applies to {\em any} $B$-module algebra to yield a formal
deformation. Known examples include formulas based on universal enveloping
algebras of Lie algebras (see examples and references in 
\cite{giaquinto-zhang98})
and a formula based on a small noncocommutative bialgebra
\cite{caldararu-giaquinto-witherspoon}.
In Section 3 we generalize the formula in 
\cite{caldararu-giaquinto-witherspoon}.
Our universal deformation formula is based on a bialgebra generated by
skew-primitive and group-like elements, and depends on a parameter $q$. 
The bialgebra and formula
were discovered in the generic case by
Giaquinto and Zhang, but were not published \cite{GZ}.
We state their formula and modify it to include the case where 
$q$ is a root of unity (Theorem \ref{thm:main}).
The case $q=-1$ is \cite[Lemma 6.2]{caldararu-giaquinto-witherspoon},
and as a first new example,
we show that the smallest Taft algebra may be deformed by applying
this earlier formula (Example \ref{ex:cibils}).
This is one of a series of algebras defined by quiver and relations
whose deformations were given by Cibils \cite{cibils90,cibils91}, and
we recover one of his deformations.
We end Section 3 with more details on the structure of the bialgebras
involved, briefly reviewing Kharchenko's construction of a Hopf algebra of
automorphisms and skew derivations of an algebra 
\cite{kharchenko91,kharchenko01}.

We apply our universal deformation formula in Section 4 to deform 
some crossed products of polynomial algebras with groups 
(Corollary \ref{cor:fd}), generalizing 
\cite[Example~6.3]{caldararu-giaquinto-witherspoon} in which the group
was $\Z/2\Z\times \Z/2\Z$.
We use Hochschild cohomology to prove that the resulting deformations
are nontrivial, by showing that their associated infinitesimals are
not coboundaries.
We apply the results of \cite{witherspoon} to show how our
deformation formula leads, in special cases, to (twisted) graded 
Hecke algebras (Example \ref{ex:tgha}).
Due to restrictive hypotheses, our formula does not give rise to very 
many of the examples of (twisted) graded Hecke algebras in
\cite{drinfeld86,etingof-ginzburg,lusztig89,ram-shepler,witherspoon}.
However it does give an infinite
series of universal deformation formulas based on noncocommutative
Hopf algebras, as well as algebras thereby undergoing nontrivial
deformation. It also proves existence of more general deformations of
certain crossed products than those that are the 
(twisted) graded Hecke algebras.

In Section 5, we give two
small examples for which our universal deformation formula nearly
provides a universal deformation of the crossed product algebra 
in the {\em other} sense
of the word {\em universal}. That is
all possible nonclassical deformations of the algebra
are parametrized by the formula, where by nonclassical we mean those
not arising from deformations of the underlying polynomial algebra itself.

For completeness, we include an appendix in which
the Hochschild cohomology of the relevant crossed product algebra is
computed.
In the case of a trivial twisting two-cocycle associated to the group,
this was done elegantly by Farinati \cite{farinati}, and by 
Ginzburg and Kaledin \cite{ginzburg-kaledin04} in
a more general geometric setting. 
Their results are easily generalized to crossed products with
twisting cocycles, however in Section 4 we need some details from a
more direct calculation.
We give an algebraic computation
similar to that in \cite{AFLS,suarez-alvarez} 
where the crossed product was taken with
a Weyl algebra instead of a polynomial algebra (and the group is symplectic).
We thank R.-O.\ Buchweitz for explaining this computation of Hochschild
cohomology to us.

We also owe many thanks to A.\ Giaquinto and J.\ Zhang for sharing their
unpublished universal 
deformation formula with us, and especially to A.\ Giaquinto 
from whom we learned algebraic deformation theory.
We thank C.\ Cibils and J.\ Stasheff for many comments on earlier versions
of this paper.

We will work over the complex numbers, although the definitions make
sense more generally. Unless otherwise indicated, $\ot=\ot_\C$.

\section{Definitions}

Let $S$ be a $\C$-algebra.
Denote by $\Aut_\C(S)$ the group 
of all $\C$-algebra automorphisms of $S$ that preserve the
multiplicative identity. Let $g,h\in\Aut_\C(S)$.
A $g,h$-{\em skew derivation}
of $S$ is a $\C$-linear function $D:S\rightarrow S$ such that 
$$
  D(rs)=D(r)g(s) + h(r) D(s)
$$
for all $r,s\in S$. If $g=h=1$ (the identity automorphism), then $D$ is simply
a derivation of $S$. 

We will be interested in skew derivations of a
crossed product algebra which we define next. 
For more details on group crossed products, see \cite{passman89}.
Let $G$ be any subgroup of $\Aut_\C(S)$.
Let $\alpha: G\times G\rightarrow \C^{\times}$ be a {\em two-cocycle}, 
that is a function satisfying
\begin{equation}\label{eqn:schur}
  \alpha(g,h)\alpha(gh,k)=\alpha(g,hk)\alpha(h,k)
\end{equation}
for all $g,h,k\in G$.  The {\em crossed product ring}
$S\#_{\alpha}G$ is $S\ot \C G$ as a vector space, with multiplication
$$
  (r\ot g)(s\ot h) = \alpha(g,h) r \cdot g( s)\ot gh
$$
for all $r,s\in S$ and $g,h\in G$.
This product is associative as $\alpha$ is a two-cocycle.
We say $\alpha$ is a {\em coboundary} if there is some
function $\beta: G\rightarrow \C^{\times}$ such that $\alpha(g,h)=
\beta(g)\beta(h)\beta(gh)^{-1}$ for all $g,h\in G$.
The set of two-cocycles modulo coboundaries forms an abelian
group under pointwise multiplication, that is $(\alpha \alpha ')(g,h)=
\alpha(g,h)\alpha'(g,h)$ for all $g,h\in G$. The crossed product algebras
$S\#_{\alpha}G$ and $S\#_{\alpha'}G$ are isomorphic if $\alpha '
=\alpha\beta$ for some coboundary $\beta$ (that is $\alpha$ and
$\alpha'$ are cohomologous).

We will abbreviate the element $r\ot g$ of $S\#_{\alpha}G$
by $r\overline{g}$.
We will assume that $\alpha$ is {\em normalized} so that 
$
  \alpha(1,g)=\alpha(g,1)=1
$
for all $g\in G$. Thus $\overline{1}$ is the multiplicative identity
of $S\#_{\alpha}G$, and it also
follows from this and (\ref{eqn:schur})
that $\alpha(g,g^{-1})=\alpha(g^{-1},g)$ for all $g\in G$.
The action of $G$ on $S$ extends to an inner action on $S\#_{\alpha}G$,
with $g(a)=\overline{g} a (\overline{g})^{-1}$ for all $g\in G$,
$a\in S\#_{\alpha}G$, where $(\overline{g})^{-1}=\alpha^{-1}(g,g^{-1})
\overline{g^{-1}}=\alpha^{-1}(g^{-1},g)\overline{g^{-1}}$.

Now let $t$ be an indeterminate. A {\em formal deformation} of a $\C$-algebra
$A$ (for example $A=S\#_{\alpha}G$) is
an associative algebra $A[[t]]=\C[[t]]\ot A$ over formal power series $\C[[t]]$ with
multiplication
\begin{equation}\label{eqn:star}
  a*b=ab+\mu_1(a\ot b)t+\mu_2(a\ot b)t^2+\cdots
\end{equation}
for all $a,b\in A$, where $ab$ denotes the product in $A$ and
the $\mu_i:A\ot A\rightarrow A$ are
$\C$-linear maps extended to be $\C[[t]]$-linear. Associativity of
$A[[t]]$ implies that $\mu_1$ is a {\em Hochschild two-cocycle}, that is
\begin{equation}\label{eqn:cocycle}
  \mu_1(a\ot b)c + \mu_1(ab\ot c)=\mu_1(a\ot bc)+a\mu_1(b\ot c)
\end{equation}
for all $a,b,c\in A$, 
as well as further conditions on the $\mu_i$, $i\geq 1$. Thus a Hochschild
two-cocycle $\mu_1$ is the first step towards a formal deformation,
and it is called the {\em infinitesimal} of the deformation.
In general it is difficult to determine whether a given
$\mu_1$ lifts to a formal deformation of $A$. 
Hochschild cohomology is defined in the appendix; for more details
on Hochschild cohomology and
deformations of algebras, see \cite{gerstenhaber-schack92}.

One way in which to obtain a formal deformation of $A$ is through
the action of a bialgebra on $A$. A {\em bialgebra} over $\C$ 
is an associative
$\C$-algebra $B$ with $\C$-algebra maps $\Delta:B\rightarrow B\ot B$
and $\varepsilon:B\rightarrow \C$ such that 
$(\Delta\otimes\id)\circ\Delta = \Delta\circ (\id\otimes\Delta)$ and
$(\varepsilon\otimes\id)\circ\Delta = \id = (\id\otimes\varepsilon)\circ
\Delta$. We will use the standard notation 
$\Delta(b)=\sum b_1\ot b_2$ for all $b\in B$, where
the subscripts are merely place-holders. A {\em Hopf algebra} is a bialgebra
$H$ with a $\C$-linear map $S:H\rightarrow H$ such that 
$\sum (Sh_1)h_2=\varepsilon(h)=\sum h_1(Sh_2)$ for all $h\in H$.
For details on bialgebras and Hopf algebras, see \cite{montgomery93}.

A {\em universal deformation formula} (or {\em UDF})
based on a bialgebra $B$ is an element $F\in (B\ot B)
[[t]]$ of the form $F=1\ot 1 + tF_1 + t^2 F_2+\cdots$ with
each $F_i\in B\ot B$, satisfying
\begin{equation}\label{eqn:udfdefn}
  (\varepsilon\ot \id)(F)=1\ot 1=(\id\ot\varepsilon)(F)
\end{equation}
$$\mbox{and } \ [(\Delta\ot \id)(F)] (F\ot 1) = [(\id\ot \Delta)(F)]
  (1\ot F),
$$
where $\id$ denotes the identity map.
These relations (\ref{eqn:udfdefn}) are similar to some of the defining
relations for an $R$-matrix, and in fact both the inverse $R^{-1}$ and the
transpose $R_{21}$ of an $R$-matrix satisfy (\ref{eqn:udfdefn})
(see \cite{benkart-witherspoon,giaquinto-zhang98}).

Suppose $A$ is a $\C$-algebra which is also a left $B$-module.
Then $A$ is a left $B$-{\em module algebra} if
\begin{equation}\label{eqn:modalg}
  h\cdot (ab) = \sum (h_1\cdot a)(h_2\cdot b)
  \ \ \mbox{ and } \ h\cdot 1 = \varepsilon (h) 1
\end{equation}
for all $a,b\in A$ and $h\in B$.
This may be extended to a $\C[[t]]$-linear action of $B[[t]]$
by extending the scalars for $A$ to $\C[[t]]$. 
Let $m:A\ot A\rightarrow A$ denote multiplication in $A$, extended
to be $\C[[t]]$-linear. The following proposition combines
Theorem 1.3 and Definition 1.13 of \cite{giaquinto-zhang98}, and
we sketch a proof for completeness.

\begin{prop}[Giaquinto-Zhang]
\label{prop:gz}
Let $B$ be a bialgebra, $A$ a left $B$-module algebra, and $F$
a universal deformation formula based on $B$. There is a formal
deformation of $A$ given by
$
  a*b = (m\circ F)(a\ot b)
$
for all $a,b\in A$.
\end{prop}

\begin{proof}
The format of $F$ as a power series in $t$
implies that $a*b=(m\circ F)(a\ot b)$ takes
the form (\ref{eqn:star}). Associativity of $*$ follows from the second
relation in (\ref{eqn:udfdefn}) and the first relation in (\ref{eqn:modalg}):
\begin{eqnarray*}
 m\circ F\circ (m\ot \id)\circ (F\ot 1) &=& m\circ (m\ot \id)\circ
  [(\Delta\ot\id)(F)]\circ (F\ot 1)\\
  &=& m\circ (\id\ot m)\circ [(\id\ot\Delta)(F)]\circ (1\ot F)\\
  &=& m\circ F\circ (\id\ot m)\circ (1\ot F)
\end{eqnarray*}
as functions from $A[[t]]\ot_{\C[[t]]}A[[t]]\ot_{\C[[t]]} A[[t]]$ to
$A[[t]]$. Note that the first relation in (\ref{eqn:udfdefn}) and the
second relation in (\ref{eqn:modalg}) imply that $1_A$ remains the
multiplicative identity under $*$.
\end{proof}

We will need the following notation. Let $q\in \C^{\times}$.
For every
integer $i\geq 1$, let 
$
  (i)_q = 1 + q + q^2+\cdots q^{i-1},
$
and set $(0)_q=0$. Let $(i)_q !=(i)_q(i-1)_q\cdots (1)_q$ 
and $(0)_q!=1$.
The {\em $q$-binomial coefficients} are
$$
  \binom{k}{i}_q =\frac{(k)_q!}{(i)_q!(k-i)_q!}
$$
for any two integers $k\geq i\geq 0$.
The well-known $q$-binomial formula states that
\begin{equation}\label{eqn:qbinom}
  (y+z)^k=\sum_{i=0}^k\binom{k}{i}_q y^iz^{k-i}
\end{equation}
in any $\C$-algebra in which $y,z$ are elements such that $zy=qyz$.

Let $S=\C[x_1,\ldots, x_n]$ and $q\in\C^{\times}$. 
As in \cite[Exer.\ IV.9.4]{kassel95},
define the linear maps $q$-{\em differentiation}
$\partial_{i,q}: S\rightarrow S$ by
\begin{equation}\label{eqn:partial}
  \partial_{i,q} (x_1^{k_1}\cdots
  x_n^{k_n}) = (k_i)_q x_1^{k_1}\cdots x_{i-1}^{k_{i-1}}x_i^{k_i-1}x_{i+1}
  ^{k_{i+1}}\cdots x_n^{k_n}.
\end{equation}
If $q=1$, these are the usual partial differentiation operators.
If $q$ is an $\ell$th root of unity, then $\partial_{i,q}^{\ell}=0$
as $(k)_q=0$ whenever $k$ is a multiple of $\ell$.
In general, $\partial_{i,q}$ is a skew derivation on $S$, specifically
$$
  \partial_{i,q}(rs)=\partial_{i,q}(r)\tau_{i,q}(s)+r\partial_{i,q}(s)
$$
for all $r,s\in S$,
where $\tau_{i,q}$ is the automorphism of $S$ defined by $\tau_{i,q}
(x_1^{k_1}\cdots x_n^{k_n})=q^{k_i}x_1^{k_1}\cdots x_n^{k_n}$.
Under some conditions, these skew derivations may be extended to a
crossed product of $S$ with a group of linear automorphisms, as we
will see.

\section{A universal deformation formula}

Let $q\in \C^{\times}$ and let $H$ be the algebra generated by $D_1,D_2,
\sigma ^{\pm 1}$, subject to the relations
$$
  D_1D_2  = D_2D_1, \ \ \ q\sigma D_i
   = D_i\sigma  \ (i=1,2) , \ \ \ \sigma\sigma^{-1}
   =  1 =  \sigma^{-1}\sigma .
$$
It is straightforward to check that $H$ is a Hopf algebra with
$$
\begin{array}{rclrclrcl}
\Delta(D_1)&=& D_1\ot \sigma + 1\ot D_1, \ & \ \varepsilon(D_1)&=&0,
  \ & \ S(D_1)&=&-D_1\sigma^{-1},\\ 
\Delta (D_2)&=& D_2\ot 1 + \sigma \ot D_2, & \varepsilon(D_2)&=&0,
  & S(D_2)&=&-\sigma^{-1}D_2,\\ 
\Delta(\sigma)&=&\sigma\ot\sigma, &\varepsilon(\sigma)&=&1,
  & S(\sigma)&=&\sigma^{-1}.
\end{array}
$$
If $q$ is a primitive $\ell$th root of unity ($\ell \geq 2$),
the ideal $I$ generated by $D_1^{\ell}$ and $D_2^{\ell}$
is a Hopf ideal, that is $\Delta(I)\subseteq I\ot H+H\ot I$, $\varepsilon
(I)=0$ and $S(I)\subseteq I$. Checking the condition on $\Delta$ involves
the $q$-binomial formula (\ref{eqn:qbinom}) and the observation that
$\displaystyle{\binom{\ell}{i}}_q=0$ whenever $1\leq i\leq \ell -1$.
Thus the quotient $H/I$ is also a Hopf algebra. 
Let 
\begin{equation}\label{eqn:hq}
  H_q = \left\{\begin{array}{rl}
   H/I, & \mbox{ if }q\mbox{ is a primitive }\ell\mbox{th root of 
  unity }(\ell\geq 2)\\
  H, & \mbox{ if }q=1\mbox{ or is not a root of unity.}
\end{array}\right.
\end{equation}

We will need the following lemma to obtain a universal deformation 
formula based on $H_q$. If $q=1$ or is not a root of unity, we define
the $q$-exponential function by
$$
  \exp _q(y) = \sum_{i=0}^{\infty} \frac{1}{(i)_q !} y^i
$$
for any element $y$ of a $\C$-algebra in which this sum is defined.
In the proof of Theorem \ref{thm:main} below,
the $\C$-algebra will be $(H_q\ot H_q\ot H_q)[[t]]$.
If $q\neq 1$ is a root of unity, this formula makes no sense as some
denominators will be zero. We modify the formula as follows in this case.
Suppose $q$ is a primitive $\ell$th root of unity for $\ell\geq 2$.
Then we define
$$
  \exp_q(y) =\sum_{i=0}^{\ell -1} \frac{1}{(i)_q !} y^i
$$
for any element $y$ of a $\C$-algebra.

\begin{lem}\label{lem:expq}
Suppose $\ell\geq 2$, $q$ is a primitive $\ell$th root of unity, 
and $y,z$ are elements of a $\C$-algebra such that $zy=qyz$ and 
$y^iz^{\ell -i}=0$
for $0\leq i\leq \ell$. Then
$$
  \exp _q(y+z) = \exp _q(y)\exp_q(z).
$$
\end{lem}

\begin{proof}
By the assumed relations and the $q$-binomial formula (\ref{eqn:qbinom}), 
each side of the desired equation may be written in the form
$$
  \sum_{i,j=0}^{\ell -1} \frac{1}{(i)_q ! (j)_q !} y^iz^j,
$$
and thus they are equal.
\end{proof}

We note that if $q$ is not a primitive root of 1, it is a standard
result that $\exp_q(y+z) = \exp_q(y)\exp_q(z)$ whenever
$zy=qyz$ and the relevant sums are defined (see for example 
\cite[Prop.\ IV.2.4]{kassel95}).
In the root of unity case, the additional hypothesis stated in the above
lemma is required.

\begin{thm}\label{thm:main}
Let $q\in\C^{\times}$ and let $H_q$ be the Hopf algebra defined in
(\ref{eqn:hq}).
Then $\exp_q(tD_1\ot D_2)$ is a universal deformation formula based on
$H_q$.
\end{thm}

\begin{proof}
In case $q$ is not a root of unity, this is an unpublished result of
Giaquinto and Zhang \cite{GZ}. 
Their proof may be adapted to the case $q$ is a
root of unity by using Lemma \ref{lem:expq} as follows.
(The proof in case $q$ is not a root of unity is essentially the same.)
Note that the hypotheses of Lemma \ref{lem:expq} hold for the pairs
$y=tD_1\ot \sigma\ot D_2$, $z=t\ot D_1\ot D_2$ and $y=tD_1\ot\sigma\ot D_2$,
$z=tD_1\ot D_2\ot 1$ as $D_2^{\ell}=0$ and $D_1^{\ell}=0$, respectively.
As $\Delta$ is an algebra homomorphism and $D_1$ commutes with $D_2$, 
we thus have
\begin{eqnarray*}
  &&\hspace{-.6in}(\Delta\ot \id)(\exp_q(tD_1\ot D_2))[\exp_q (tD_1\ot D_2)\ot 1] \\
  &=&\exp_q((\Delta\ot \id)(tD_1\ot D_2))\exp_q(tD_1\ot D_2\ot 1)\\
  &=&\exp_q(tD_1\ot\sigma\ot D_2 + t\ot D_1\ot D_2)\exp_q(tD_1\ot D_2\ot 1)\\
  &=& \exp_q(tD_1\ot \sigma\ot D_2)\exp_q(t\ot D_1\ot D_2)\exp_q(tD_1\ot D_2\ot 1)\\
  &=& \exp_q(tD_1\ot \sigma\ot D_2)\exp_q(tD_1\ot D_2\ot 1)\exp_q(t\ot D_1\ot D_2)\\
  &=& \exp_q(tD_1\ot \sigma \ot D_2 + tD_1\ot D_2\ot 1)\exp_q(t\ot D_1\ot D_2)\\
  &=& \exp_q((\id\ot \Delta)(tD_1\ot D_2))(1\ot\exp_q(tD_1\ot D_2))\\
  &=& (\id\ot \Delta)(\exp_q(tD_1\ot D_2))(1\ot \exp_q(tD_1\ot D_2)).
\end{eqnarray*}
The remaining relation in (\ref{eqn:udfdefn}) holds as $\varepsilon(D_1)=
\varepsilon(D_2)=0$. Thus $\exp_q(tD_1\ot D_2)$ is a universal
deformation formula.
\end{proof}

By Proposition \ref{prop:gz} and Theorem \ref{thm:main}, we need only
find an $H_q$-module algebra $A$, and $m\circ \exp _q (tD_1\ot D_2)$ 
will provide a formal deformation of $A$.
Our first such example is next; a large family of examples is given in
Section 4.

\begin{example}\label{ex:cibils}{\em (A Taft algebra.)
Let $A$ be the algebra defined by generators and relations as follows,
where the indices are read modulo $2$:
\begin{eqnarray*}
  A=\C\langle s_0,s_1,\gamma_0,\gamma_1 &\mid & s_0+s_1=1, \
  s_i^2=s_i, \ s_is_{i+1}=0, \gamma_i^2=0, \ \gamma_i\gamma_{i+1}=0,\\
  && \hspace{.6in} s_i\gamma_i =0, \ s_{i+1}\gamma_i=\gamma_i,
  \ \gamma_is_i=\gamma_i, \ \gamma_is_{i+1}=0\rangle .
\end{eqnarray*}
This is an algebra defined by a quiver and relations
as in \cite[Thm.\ 5.1(b)]{cibils90}; 
the quiver is Gabriel's quiver consisting of two arrows in opposite
directions between two vertices.
The algebra $A$ is isomorphic to 
$\C \langle x,g\mid gx=-xg, \ x^2=0, \ g^2=1\rangle $ 
via the map $x\mapsto \gamma_0-\gamma_1$, $g\mapsto s_0-s_1$.
It has the structure of a Hopf algebra first discovered by Sweedler,
and is one of a series of Hopf algebras constructed by Taft
\cite[Example 1.5.6]{montgomery93}.
In \cite{cibils90,cibils91}, 
Cibils gave deformations of more general classes of
algebras defined by quivers and relations, and in this special case one of his
deformations may be obtained by applying a universal deformation formula.
Specifically, let $q=-1$ and
$$
  H_{-1}=\C\langle D_1,D_2,\sigma^{\pm 1}\mid D_1D_{2}=D_{2}D_1, \ -\sigma
  D_i=D_i\sigma, \ D_i^2=0, \ \sigma\sigma^{-1}=1=\sigma^{-1}\sigma\rangle
$$
as above.
Define
$$
 \sigma(\gamma_i)=-\gamma_{i+1}, \ \sigma(s_i)=s_{i+1}, 
$$
$$
  D_1(\gamma_i)=s_{i+1}, \ D_2(\gamma_i)=s_i, \ D_i(s_j)=0.
$$
It may be checked that the relations of $H_{-1}$ 
are preserved on the generators of $A$, making the vector space
$V=\Span _{\C}\{s_0,s_1,\gamma_0,\gamma_1\}$ into an $H_{-1}$-module.
Therefore the tensor algebra $T(V)$ is an $H_{-1}$-module algebra,
where the action of $H_{-1}$ is extended to $T(V)$ by (\ref{eqn:modalg}).
As $A$ is a quotient of $T(V)$, it remains to check that the relations
of $A$ are preserved by the generators of $H$, a straightforward 
computation.
(In fact, $A$ is also a module algebra for the finite dimensional 
quotient $H_{-1}/(\sigma^2 - 1)$.)
By Proposition \ref{prop:gz} and 
Theorem \ref{thm:main}, $\exp_{-1}(tD_1\ot D_2)
=1+tD_1\ot D_2$ yields a formal deformation of $A$. The deformation is
\begin{eqnarray*}
  A_t = \C\langle s_0,s_1,\gamma_0,\gamma_1 &\mid &
   s_0+s_1=1, \ s_i^2=s_i, \ s_is_{i+1}=0, \ \gamma_i^2=0, \ s_i\gamma_i=0,\\
   && \hspace{.1in} s_{i+1}\gamma_i=\gamma_i, \ \gamma_is_i=\gamma_i,
   \gamma_is_{i+1}=0, \ \gamma_i\gamma_{i+1}=ts_{i+1}\rangle ,
\end{eqnarray*}
which is precisely that given in \cite[Thm.\ 5.1(b)]{cibils90}.
This deformation is nontrivial since if we specialize to $t\neq 0$,
$A_t$ is isomorphic to the $2\times 2$ matrix algebra, and thus is not 
isomorphic to $A$.
We do not know whether any of Cibils' other deformations are given by
universal deformation formulas.
}
\end{example}

More generally, suppose that $A$ is any $H_q$-module algebra.
Due to (\ref{eqn:modalg}) and the nature of the coproducts of $D_1,D_2$, 
the following general lemma implies that $\mu_1 = m\circ (D_1\ot D_2)$
is a Hochschild two-cocycle on $A$, that is it satisfies (\ref{eqn:cocycle}).
This generalizes the well-known fact that the cup product of derivations
is a Hochschild two-cocycle.
The lemma is proved by direct computation, with no assumption made on
the relations among $D_1,D_2,\sigma$.
If the relations of $H_q$ do hold however, then Theorem \ref{thm:main} gives
an alternative proof that $\mu_1=m\circ (D_1\ot D_2)$ is
a Hochschild two-cocycle.

\begin{lem}\label{lem:htc}
Let $A$ be an algebra over a field, with multiplication
$m: A\otimes A\rightarrow A$. Let $\sigma$ be an automorphism of $A$,
$D_1$ a $\sigma, 1$-skew derivation and $D_2$ a $1,\sigma$-skew
derivation of $A$. Then $\mu_1=m\circ (D_1\otimes D_2)$ is a Hochschild
two-cocycle.
\end{lem}

We end this section with
a construction due to Kharchenko \cite[\S6.5.5]{kharchenko91}
of a Hopf algebra of automorphisms and skew derivations of an algebra $A$.
The Hopf algebras $H_q$ are related to some of 
Kharchenko's Hopf algebras,
and it may be useful to consider his general
construction in questions regarding deformations of algebras.

Let $K$ be a subgroup of $\Aut_{\C}(A)$.
For each $k\in K$, let $L_k$
be a vector subspace of $\End_{\C}(A)$
consisting of $1,k$-skew derivations of $A$. 
Let $L=\oplus_{k\in K}L_k$, and assume $K$ acts on $L$ in such a way that
$m(L_k)=L_{mkm^{-1}}$ for all $k,m\in K$.
Thus $K$ acts by automorphisms on the tensor algebra $T(L)$, 
and we let $H=T(L) \# K$. The coproducts 
$$
  \Delta(k)=k\ot k \ \ \mbox{ and } \ \
  \Delta (D) =D\ot 1 + k\ot D
$$
for all $k\in K$ and $D\in L_k$ extend, by requiring $\Delta$ to be an algebra
homomorphism, to a coproduct $\Delta$ on $H$. Similarly, the counit
$\varepsilon$ and antipode $S$ defined as follows on generators
extend to $H$: $\varepsilon(k)=1$, $\varepsilon(D)=0$,
$S(k)=k^{-1}$ and $S(D)=-k^{-1}D$ for all $k\in K$ and $D\in L_k$.
Thus $H$ is a Hopf algebra.

We obtain the Hopf algebras $H_q$ by this construction in the following way:
If $A$ is an $H_q$-module algebra,
let $K$ be the group generated by the action of $\sigma$ on $A$,
$L_{\sigma}=\Span_{\C}\{D_2\}$, $L_{\sigma^{-1}}=\Span_{\C}
\{D_1 \sigma^{-1}\}$, and $L_{\tau}=\{0\}$ if $t\neq \sigma^{\pm 1}$.
If $\sigma$ has infinite order as an automorphism of $A$, 
then $H_q$ is a quotient of Kharchenko's Hopf algebra defined by this data. 
Otherwise we must take a quotient of $H_q$, in which $\sigma$
has the correct order, to obtain a quotient of Kharchenko's Hopf algebra.
See \cite{kharchenko01} for further details on this construction.

\section{Deformations of group crossed products}

In this section we give a large family of group crossed products
to which the formula
of Theorem \ref{thm:main} applies to yield nontrivial deformations.

Let $G$ be a group with a representation on a $\C$-vector space $V$
of dimension $n$, so that $G$ acts by automorphisms on the symmetric
algebra $S(V)$. We will identify $S(V)$ with polynomials in the
variables $x_1,\ldots, x_n$.
In this section, we will be interested in formal deformations of
a crossed product $S(V)\#_{\alpha}G$ for which the infinitesimal $\mu_1$
satisfies $\mu_1(V\ot V)\subset S(V)\overline{g}$ for some $g\in G$.
Not all elements $g\in G$ correspond to such noncoboundary 
infinitesimals $\mu_1$.
In case $G$ is finite, examination
of Hochschild cohomology (see Corollary \ref{cor:additive} and 
subsequent comments) shows that we may assume such an element $g$ 
has determinant 1 as an operator on $V$,
and $\codim(V^g)=0$ or $2$, where $V^g=\{v\in V\mid g(v)=v\}$, 
the subspace of $V$ invariant under $g$.
In this section, we will make this assumption,
and in addition will assume that $g$ is central in $G$. 
Again if the order of $g$ is finite, $g$ acts diagonally
with respect to some basis of $V$, and without loss of generality this is 
$x_1,\ldots,x_n$. Specifically, we will assume that
\begin{equation}\label{eqn:detg}
  g(x_1)=qx_1, \ g(x_2)=q^{-1}x_2, \ g(x_3)=x_3, \cdots , \ g(x_n)=x_n
\end{equation}
for some $q\in \C^{\times}$.
In order to include some infinite groups, we will not assume that $q$ is a root
of unity. To obtain explicit formulas, we will further need to make
a more restrictive assumption: 
\begin{equation}\label{eqn:preserves}
G \mbox{ preserves the subspaces } \C x_1, \C x_2 \mbox{ of } V. 
\end{equation}
If $q\neq \pm 1$, this is automatically
the case by the assumed centrality of $g$. 
Under the assumption (\ref{eqn:preserves}),
we may abuse notation and define the functions
$x_i: G\rightarrow \C^{\times}$ ($i=1,2$) by
$$
  h(x_i)=x_i(h)x_i
$$
for each $h\in G$.
Let $D_1,D_2$ and $\sigma:S(V)\#_{\alpha} G\rightarrow S(V)\#_{\alpha} G$
be the linear functions defined on a basis 
$\{x_1^{k_1}\cdots x_n^{k_n}\overline{h}\mid k_i\in \Z^{\geq 0}, \ h\in G\}$
as follows:
\begin{eqnarray}
   D_1(x_1^{k_1}\cdots x_n^{k_n}\overline{h})& =& x_1(h^{-1})\partial_{1,q}
  (x_1^{k_1}\cdots x_n^{k_n})\overline{h},
    \label{eqn:d1}\\
  D_2(x_1^{k_1}\cdots x_n^{k_n}\overline{h})& =& q^{k_1} \partial_{2,q^{-1}}
  (x_1^{k_1}\cdots x_n^{k_n}) 
  s\overline{g}\cdot \overline{h},
    \label{eqn:d2}\\
  \sigma(x_1^{k_1}\cdots x_n^{k_n}\overline{h})& =& x_1(h^{-1}) q^{k_1}x_1^{k_1}\cdots 
  x_n^{k_n} \overline{h},
    \label{eqn:beta}
\end{eqnarray}
where $\partial _{1,q}$, $\partial_{2,q^{-1}}$ are defined in
(\ref{eqn:partial}) and $s\in \C[x_3,\ldots,x_n]$ satisfies
\begin{equation}\label{eqn:s}
  h(s) = x_1(h)x_2(h)\alpha(g,h)\alpha^{-1}(h,g)s
\end{equation}
for all $h\in G$, that is $s$ is a semi-invariant of $G$.
(Our condition on the polynomial $s$ is
informed by knowledge of Hochschild cohomology; see Corollary 
\ref{cor:additive} and the computation (\ref{eqn:sinvariant}) below.)
Calculations using (\ref{eqn:detg})--(\ref{eqn:s}) and centrality of $g$
show that $\sigma$ is an automorphism and $D_1,D_2$ are skew derivations
with respect to $\sigma$, specifically 
$$
  D_1(ab)= D_1(a)\sigma(b) + aD_1(b) \ \ \ \mbox{ and } \ \ \
  D_2(ab) = D_2(a)b + \sigma(a)D_2(b)
$$
for all $a,b \in S(V)\#_{\alpha} G$. A direct calculation shows that 
$\mu_1=m\circ (D_1\otimes D_2)$ is a Hochschild two-cocycle on $S(V)\#_{\alpha} G$,
that is $\mu_1$ satisfies (\ref{eqn:cocycle}).
This is also a consequence of Lemma \ref{lem:htc}, or of 
Theorem \ref{thm:main} in combination with Theorem \ref{thm:modalg} below.

Taking $D_1,D_2$ to be the skew derivations defined in (\ref{eqn:d1}),
(\ref{eqn:d2}), the corresponding Hochschild two-cocycle $\mu_1$
takes $V\otimes V$ to 
$S(V)\overline{g}$, the $g$-component of $S(V)\#_{\alpha}G$. 
If $g$ were {\em not} central in $G$, an associated Hochschild two-cocycle
would necessarily involve all components of $S(V)\#_{\alpha}G$ 
corresponding to the elements of the conjugacy class of $g$
(see Corollary \ref{cor:additive}).
We do not know if the explicit formulas of this section can be generalized
to noncentral $g$.

Let $H_q$ be the Hopf algebra defined in (\ref{eqn:hq}).

\begin{thm}\label{thm:modalg}
Let $g$ be a central element of $G$ such that
(\ref{eqn:detg}) and (\ref{eqn:preserves}) hold. Then
$S(V)\#_{\alpha} G$ is an $H_q$-module algebra under the action defined in
(\ref{eqn:d1})--(\ref{eqn:beta}).
\end{thm}

\begin{proof}
The relations among the generators in $H_q$ may be checked to be 
preserved under the action, so that $S(V)\#_{\alpha} G$ is an $H_q$-module.
In particular, in case $q$ is a primitive $\ell$th root of unity,
$D_1^{\ell}=0=D_2^{\ell}$ as $(k)_q=0$ whenever $k$ is a multiple
of $\ell$.
As stated earlier, $\sigma$ is an automorphism of $S(V)$ and $D_1$ and $D_2$
are skew derivations. Clearly
$D_1(1)=0=D_2(1)$ as $(0)_q=0$. 
Therefore (\ref{eqn:modalg}) holds, so
$S(V)\#_{\alpha} G$ is an $H_q$-module algebra.
\end{proof}

Combining Proposition \ref{prop:gz} and
Theorems \ref{thm:main} and \ref{thm:modalg},
we now have the following corollary.
In case $G$ is finite, the deformations in the corollary are shown to
be nontrivial in the remainder of this section.
We expect that the same is also true in case $G$ is infinite.

\begin{cor}\label{cor:fd}
Let $g$ be a central element of $G$ such that (\ref{eqn:detg}) and
(\ref{eqn:preserves}) hold. Then
$\exp_q(tD_1\ot D_2)$ yields a formal deformation of $S(V)\#_{\alpha} G$.
\end{cor}

We point out that if $g=1$ then $\exp_q(tD_1\ot D_2)$ restricts to a 
classical formula on $S(V)$, namely
$$
  \sum_{i=0}^{\infty} \frac{t^i}{i!} \left(\frac{\partial}{\partial
   x_1}\right)^i\ot \left(s\frac{\partial}{\partial x_2}\right)^i,
$$
where $s\in \C[x_3,\ldots,x_n]$ satisfies $h(s)=x_1(h)x_2(h)s$ for all
$h\in G$. Taking $G$ to be the identity group, $n=2$, and $s=1$,
this formula applied to $\C[x_1,x_2]$ yields the Weyl algebra on
two generators.

In the special case 
$G=\Z/2\Z\times \Z/2\Z$ and $V=\C^3$ with a particular diagonal
action of $G$,
the formula and deformation of Corollary \ref{cor:fd} were obtained in 
\cite[\S6]{caldararu-giaquinto-witherspoon}. 
The deformations in that case are nontrivial since
their corresponding Hochschild two-cocycles are not coboundaries, a
consequence of the computations in \cite{caldararu-giaquinto-witherspoon}.
Similarly, we now show that the same is true in the more general setting
of a finite group $G$, based on a computation of the
Hochschild cohomology of $S(V)\#_{\alpha}G$.
The Hochschild cohomology was computed by Farinati,
Ginzburg and Kaledin in the case 
$\alpha =1$ \cite{farinati,ginzburg-kaledin04}. The addition of a nontrivial
cocycle $\alpha$ poses no difficulties, however we need to use
some of the details from an explicit algebraic computation.
These we provide in the appendix. 

There is a chain map from the bar complex (\ref{eqn:bar}) for $A=S(V)$
to the Koszul complex $K(\{x_i\ot 1 - 1\ot x_i\}_{i=1}^n)$,
$$
\begin{array}{ccccccccccc}
\cdots \! &\rightarrow & S(V)^{\ot 4} & \stackrel{\delta_2}{\longrightarrow}
  & S(V)^{\ot 3} & \stackrel{\delta_1}{\longrightarrow} & S(V)^e &
  \stackrel{m}{\longrightarrow} & S(V) & \rightarrow & \! 0\\
 & & \hspace{.1in}\downarrow \psi_2 & & \hspace{.1in}\downarrow \psi_1 & &
   \parallel &&\parallel &&\\
 \cdots \! &\! \rightarrow \! & \Wedge^2(V)\ot S(V)^e \! & \! 
  \stackrel{d_2}{\longrightarrow}\!
  & \!\Wedge^1(V)\ot S(V)^e \! &\! \stackrel{d_1}{\longrightarrow} 
  \! & \! S(V)^e &\stackrel
  {m}{\longrightarrow}\! & \! S(V) \! &\! \rightarrow & \! 0 
\end{array}
$$
We will need an explicit formula for $\psi_2$ in particular.
A straightforward computation shows that the following formulas work
(cf.\ \cite{caldararu-giaquinto-witherspoon}, in which slightly
different formulas are given in the case $n=3$):
\begin{equation}\label{eqn:psi1}
 \psi_1(1\ot x_1^{k_1}\cdots x_n^{k_n}\ot 1) = \sum_{i=1}^n \sum_{p=1}
   ^{k_i}e_i\ot x_i^{k_i-p}x_{i+1}^{k_{i+1}}\cdots x_n^{k_n}\ot 
  x_1^{k_1}\cdots x_{i-1}^{k_{i-1}} x_i^{p-1},
\end{equation}
\begin{equation}\label{eqn:psi2}
\psi_2(1\ot x_1^{k_1}\cdots x_n^{k_n}\ot x_1^{m_1}\cdots x_n^{m_n}\ot 1)=
\hspace{4in}
\end{equation}
$$\sum_{1\leq i<j\leq n}\sum_{r=1}^{m_j}\sum_{p=1}^{k_i}  e_i\wedge e_j\ot
  x_i^{k_i-p} x_{i+1}^{k_{i+1}}\cdots x_{j-1}^{k_{j-1}}x_j^{k_j+m_j-r}x_{j+1}
  ^{k_{j+1}+m_{j+1}}\cdots x_n^{k_n+m_n}\ot 
$$

\vspace{-.15in}

$$\hspace{2in}x_1^{k_1+m_1}\cdots x_{i-1}^
 {k_{i-1}+m_{i-1}} x_i^{m_i+p-1} x_{i+1}^{m_{i+1}}\cdots x_{j-1}^{m_{j-1}}
  x_j^{r-1}.
$$

Now assume $G$ is a finite group acting on $V$, and $g$ is a central
element of $G$ satisfying (\ref{eqn:detg}) and (\ref{eqn:preserves}) 
where $q$ is a primitive $\ell$th root of unity, $\ell \geq 2$.
Under these assumptions, by Proposition \ref{thm:additive} and
Corollary \ref{cor:additive}, $\HH^2(S(V)\#_{\alpha}G)$ contains
as the $g$-component
$$
  \HH^2(S(V),S(V)\overline{g})^G\cong 
    (\det(\Span_{\C}\{x_1,x_2\}^*)\ot \C[x_3,\ldots x_n]\overline{g})^G.
$$
Let $s\in \C[x_3,\ldots,x_n]-\{0\}$ satisfy (\ref{eqn:s}),
that is $h(s)=x_1(h)x_2(h)\alpha(g,h)\alpha^{-1}(h,g)s$ for all
$h\in G$. Identify the dual function $ (e_1\wedge e_2)^*$ with a
basis of the one-dimensional space $\det(\Span_{\C}\{x_1,x_2\}^*)$, where
the notation $e_i$ comes from the Koszul complex and is defined in
the appendix.
We first
claim that $ (e_1\wedge e_2)^* \ot s\overline{g}$ corresponds to a nonzero 
element of $\HH^2(S(V)\#_{\alpha}G)$ under the above isomorphism.
We need only show that $(e_1\wedge e_2)^*\ot s\overline{g}$ is invariant
under the action of $G$. Let $h\in G$. Then
\begin{equation}\label{eqn:sinvariant}
  \hspace{-.118in}
   h ((e_1\wedge e_2)^*\ot s\overline{g}) = x_1(h^{-1})x_2(h^{-1})
  (e_1\wedge e_2)^*\ot (h(s))\overline{h}\cdot \overline{g}\cdot
  (\overline{h})^{-1}
\end{equation}

\vspace{-.2665in}

\begin{eqnarray*}
  \hspace{1.36in}&=& 
  \alpha(g,h)\alpha^{-1}(h,g) (e_1\wedge e_2)^* \ot s\overline{h}
  \cdot \overline{g}\cdot (\overline{h})^{-1}\\
  &=& \alpha(g,h) (e_1\wedge e_2)^*\ot s\overline{hg} \cdot (\overline{h}
  )^{-1}\\ 
  &=& \alpha(g,h)\alpha^{-1}(h,h^{-1})\alpha(hg,h^{-1}) (e_1\wedge e_2)^*
  \ot s\overline{hgh^{-1}}\\
  &=& (e_1\wedge e_2)^* \ot s\overline{g}
\end{eqnarray*}
by an application of the two-cocycle identity (\ref{eqn:schur}) to the
triple $g,h,h^{-1}$, since $g\in C(G)$.

Next we show that the nonzero element $(e_1\wedge e_2)^*\ot s\overline{g}$
of $\HH^2(S(V)\#_{\alpha}G)$ may be identified with a Hochschild
two-cocycle $\mu_1$ of the form $m\circ (D_1\ot D_2)$ where $D_1,D_2$
are defined in (\ref{eqn:d1}), (\ref{eqn:d2}).
This will prove that the formal deformations of $S(V)\#_{\alpha}G$
given by Corollary \ref{cor:fd} are nontrivial in case $G$ is finite.
We will need \cite[Thm.\ 5.4]{caldararu-giaquinto-witherspoon}, which will
be applied to a Koszul resolution:

\begin{prop}[Caldararu-Giaquinto-Witherspoon]
Let $A=S(V)\#_{\alpha}G$.
Let $f:P_n\rightarrow A$ be a function representing an element of
$\HH^n(S(V),A)^G\cong \HH^n(A)$ expressed in terms of any $S(V)^e$-projective
resolution $P_{\DOT}$ of $S(V)$ carrying an action of $G$. 
The corresponding function 
$\widetilde{f}\in \Hom_{\C}(A^{\ot n},A)\cong
\Hom_{A^e}(A^{\ot (n+2)}, A)$ on the bar complex (\ref{eqn:bar}) is given
by
$$
 \widetilde{f} (p_1\overline{\sigma_1}\ot \cdots\ot p_n\overline{\sigma_n})
  = ((f\circ \psi_n)(1\ot p_1\ot \sigma_1(p_2)\ot \cdots\ot (\sigma_1\cdots
\sigma_{n-1})(p_n)\ot 1 )) \overline{\sigma_1}\cdots \overline{\sigma_n}.
$$
In particular, if $n=2$, we obtain the infinitesimal deformation
$\mu_1 : A\ot A\rightarrow A$,
$$
  \mu_1 (p_1\overline{\sigma_1}\ot p_2\overline{\sigma_2})=((f\circ \psi_2)
  (1\ot p_1\ot \sigma_1(p_2)\ot 1))\overline{\sigma_1}\cdot \overline{\sigma_2}.
$$
\end{prop}

As a consequence of the proposition, the element $(e_1\wedge e_2)^*\ot 
s\overline{g}$ of $\HH^2(S(V)\#_{\alpha}G)$ may be identified with the
function $\mu_1:A\ot A\rightarrow A$ where 
$\mu_1(x_1^{k_1}\cdots x_n^{k_n}\overline{h}\ot x_1^{m_1}\cdots x_n^{m_n}
\overline{k})$ is $\psi_2(1\ot x_1^{k_1}\cdots x_n^{k_n}\ot
h(x_1^{m_1}\cdots x_n^{m_n})\ot 1)$ followed by application of the
function representing $(e_1\wedge e_2)^* \ot s\overline{g}$ at the chain
level, and right multiplication by $\overline{h}\cdot \overline{k}$.
By our hypotheses, we have

\noindent
$
 \psi_2(1\ot x_1^{k_1}\cdots x_n^{k_n}\ot h(x_1^{m_1}\cdots x_n^{m_n})\ot 1)
$
$$\hspace{.5in} = x_1(h)^{m_1} x_2(h)^{m_2}\psi_2(1\ot x_1^{k_1}\cdots 
 x_n^{k_n}\ot x_1^{m_1}x_2^{m_2} h(x_3^{m_3}\cdots x_n^{m_n})).
$$
By (\ref{eqn:psi2}), the resulting coefficient of $e_1\wedge e_2$ is
$$
  x_1(h)^{m_1}x_2(h)^{m_2} \sum_{r=1}^{m_2}\sum_{p=1}^{k_1} x_1^{k_1-p} x_2
  ^{k_2+m_2-r} x_3^{k_3}\cdots x_n^{k_n}\cdot h(x_3^{m_3}\cdots x_n^{m_n})
  \ot x_1^{m_1+p-1}x_2^{r-1}.
$$
Applying $(e_1\wedge e_2)^* \ot s\overline{g}$ and multiplying by
$\overline{h}\cdot \overline{k}$, we obtain 
\begin{eqnarray*}
 && \hspace{-.3in}x_1(h)^{m_1}x_2(h)^{m_2}
  \sum_{r=1}^{m_2}\sum_{p=1}^{k_1} x_1^{k_1-p} 
 x_2^{k_2+m_2-r}x_3^{k_3}\cdots x_n^{k_n} h(x_3^{m_3}\cdots x_n^{m_n}) 
  s\overline{g} x_1^{m_1+p-1} x_2^{r-1}\overline{h}\cdot
  \overline{k}\hspace{.3in}\\
&=& \!\! x_1(h)^{m_1}x_2(h)^{m_2}\sum_{r=1}^{m_2}\sum_{p=1}^{k_1} q^{m_1+p-r}
  x_1^{k_1+m_1-1}x_2^{k_2+m_2-1}x_3^{k_3}\cdots x_n^{k_n}
  h(x_3^{m_3}\cdots x_n^{m_n}) s \overline{g}\cdot\overline{h}\cdot
  \overline{k}\\
&=& \!\! q^{m_1} x_1(h)^{m_1}x_2(h)^{m_2} (k_1)_q(m_2)_{q^{-1}} 
  x_1^{k_1+m_1-1}x_2^{k_2+m_2-1} x_3^{k_3}\cdots x_n^{k_n}
  h(x_3^{m_3}\cdots x_n^{m_n}) s\overline{g}\cdot \overline{h}
  \cdot \overline{k}.
\end{eqnarray*}
On the other hand,
\begin{eqnarray*}
&&\hspace{-.3in}(m\circ (D_1\ot D_2))(x_1^{k_1}\cdots x_n^{k_n}\overline{h}
 \ot x_1^{m_1}\cdots x_n^{m_n}\overline{k})\hspace{.5in}\\
&=& \!\! x_1(h^{-1})(k_1)_q x_1^{k_1-1} x_2^{k_2}\cdots x_n^{k_n}\overline{h}
  \cdot q^{m_1} (m_2)_{q^{-1}} x_1^{m_1}x_2^{m_2-1} x_3^{m_3}
  \cdots x_n^{m_n} s\overline{g}\cdot \overline{k}\\
&=& \!\! q^{m_1}x_1(h)^{m_1}x_2(h)^{m_2}(k_1)_q(m_2)_{q^{-1}}
  x_1^{k_1+m_1-1}x_2^{k_2+m_2-1} x_3^{k_3}\cdots x_n^{k_n}
  h(x_3^{m_3}\cdots x_n^{m_n}) s\overline{g}\cdot
  \overline{h}\cdot \overline{k}.\\
\end{eqnarray*}
Therefore $m\circ (D_1\ot D_2)$ is the Hochschild two-cocycle represented by
$(e_1\wedge e_2)^*\ot s\overline{g}$ and so is not a coboundary.
This implies that the formal deformations given by Corollary \ref{cor:fd}
are nontrivial in case $G$ is finite.

\quad

\begin{example}\label{ex:tgha}
{\em (Twisted graded Hecke algebras.)
Let $G$ be a finite subgroup of $GL(V)$,
and $g\in G$ a central element satisfying (\ref{eqn:detg}) and (\ref{eqn:preserves}).
Suppose $s=1$
satisfies (\ref{eqn:s}), that is $1=x_1(h)x_2(h)\alpha(g,h)\alpha^{-1}(h,g)$
for all $h\in G$. We may rewrite this condition as $\det(h|_{(V^g)^{\perp}})
=\alpha(h,g)\alpha^{-1}(g,h)$.
Then $\mu_1=m\circ (D_1\ot D_2)$ is a bilinear form on $S(V)\#_{\alpha}G$
of degree $-2$, where $S(V)\#_{\alpha}G$ is a graded algebra in which
elements of $V$ have degree 1 and elements of $G$
have degree 0. More generally, in the formula $\exp_q(tD_1\ot D_2)$, the
bilinear form $\mu_i=\frac{1}{(i)_q!} m \circ (D_1^i\ot D_2^i)$
has degree $-2i$.
By \cite[Thm\ 3.2]{witherspoon}, 
the resulting formal deformation of $S(V)\#_{\alpha}G$
becomes a (twisted) graded Hecke algebra \cite{witherspoon}
when the scalars are restricted
to $\C[t]$. In this case, that means the associated deformation of 
$S(V)\#_{\alpha}G$ over $\C [t]$ is isomorphic to
$$
  T(V)\#_{\alpha}G[t]/(vw-wv-a_g(v,w)t\overline{g}),
$$
the quotient by the ideal generated by all elements $vw-wv-a_g(v,w)t
\overline{g}$, for 
$v,w\in V$, where $a_g(v,w)=\mu_1(v,w)-\mu_1(w,v)$. This (twisted)
graded Hecke algebra is special in that only one such function $a_g$
is nonzero. In the next section, we give some examples for which there is
an analogous deformation with
more than one group element $g$ having $a_g$ nonzero.}
\end{example}

\section{Universal Deformations}

In this section we give examples for which some of the universal
deformation formulas from the last section, corresponding to different
group elements, may be combined into larger formulas. The first
example generalizes \cite[Lemma 6.2]{caldararu-giaquinto-witherspoon}.

\begin{example}\label{ex:1}{\em
Let $n\geq 3$ and $\ell\geq 2$ be integers,
$q$ a primitive $\ell$th root of unity,  
$G=(\Z/\ell\Z)^{n-1}$ and $V=\C ^n$. Identify $G$ with the subgroup of ${\rm SL}(V)$
generated by the diagonal matrices
\begin{eqnarray*}
 g_1&=&\diag(q,q^{-1},1,\ldots , 1), \\
  g_2&=&\diag(1,q,q^{-1},1,\ldots , 1),\\
     &\vdots &\\ 
  g_{n-1}&=&\diag(1,\ldots , 1, q,q^{-1}),
\end{eqnarray*}
with respect to a basis $x_1,\ldots,x_n$ of $V$. 
Let $g_n=g_1^{-1}\cdots g_{n-1}^{-1}=\diag (q^{-1},1,\ldots, 1,q)$.
Let $\alpha: G\times G\rightarrow \C^{\times}$ be
the following two-cocycle:
$$
  \alpha(g_1^{i_1}\cdots g_{n-1}^{i_{n-1}}, \ g_1^{j_1}\cdots
   g_{n-1}^{j_{n-1}}) = q^{-\sum_{1\leq k\leq n-2} i_k j_{k+1}}.
$$
It may be checked directly that $\alpha$ satisfies the two-cocycle condition
(\ref{eqn:schur}).
Note that $\alpha$ is not a coboundary: By their definition, two-coboundaries
for abelian groups are symmetric, but
$\alpha$ is clearly not symmetric.
In case $n=3, \ \ell =2$, $\alpha$ is cohomologous to the 
nontrivial cocycle given in
\cite[Example 3.4]{caldararu-giaquinto-witherspoon}, as in that case there is a
unique nontrivial two-cocycle up to coboundary.
Note that 
\begin{equation}\label{eqn:abc}
  \alpha(g_{i+1},g_i)= q\alpha(g_i,g_{i+1}) \ \ \ \mbox{for } 1\leq i\leq n
\end{equation}
(where $g_{n+1}=g_1$ by definition).
Direct calculations also show that
\begin{equation}\label{eqn:xyz}
  h(x_i x_{i+1}) = \frac{\alpha(h,g_i)}{\alpha(g_i,h)} x_ix_{i+1} \ \ \
   \mbox{for } 1\leq i\leq n,
\end{equation}
for all $h\in G$ (where $x_{n+1}=x_1$ by definition).

Let $H_i = H_{q}$ ($1\leq i\leq n$) be the Hopf algebra 
defined in (\ref{eqn:hq}), acting on $S(V)\#_{\alpha}G$ via
the formulas (\ref{eqn:d1})--(\ref{eqn:beta}), where we replace
$g$ by $g_i$ and $x_1,x_2$ by $x_i,x_{i+1}$. 
Applying (\ref{eqn:s}) and (\ref{eqn:xyz}), 
it may be checked
that the polynomial $s$ arising in the action of $H_i$ 
must be in $S(V^{g_i})\cap (S(V))^G = \C[x_1^{\ell},\ldots , x_{i-1}^{\ell},
x_{i+2}^{\ell},\ldots , x_n^{\ell}]$, where if $i=n$ we leave out $x_1$
and $x_n$.
Using this fact, equation (\ref{eqn:abc}), and the
identity $(i+\ell)_{q^{-1}}=(i)_{q^{-1}}$ for all integers $i$, it may
be checked directly that the images of the $H_i$ ($1\leq i\leq n$) in
$\End_{\C}(S(V)\#_{\alpha}G)$ mutually commute.
Thus there is a corresponding algebra homomorphism from the Hopf algebra 
$H_1\ot \cdots
\ot H_n$ to $\End_{\C}(S(V)\#_{\alpha}G)$, and $S(V)\#_{\alpha}G$ is
a module algebra for $H_1\ot\cdots\ot H_n$. 
A product of universal deformation formulas is again a universal deformation
formula, based on the tensor product of the bialgebras. Thus by
Theorem \ref{thm:main},
\begin{equation}\label{eqn:udf}
  \exp_q(tD_1^{g_1}\ot D_2^{g_1})\cdots\exp_q(tD_1^{g_n}\ot D_2^{g_n})
\end{equation}
(where superscripts indicate the Hopf subalgebra from which the operators
originate)
is a universal deformation formula based on $H_1\ot \cdots\ot H_n$.
By Proposition \ref{prop:gz}, this formula applies to yield a formal 
deformation of $S(V)\#_{\alpha}G$.

In case $s=1$, restricting this deformation to one over $\C [t]$
results in the twisted graded Hecke algebra
$$
  T(V)\#_{\alpha}G[t]\left/ \left(vw-wv-\sum_{i=1}^n a_{g_i}(v,w)t
  \overline{g}_i\right)\right.
$$
where $a_{g_i}(v,w)=D_1^{g_i}(v)D_2^{g_i}(w)-D_1^{g_i}(w)D_2^{g_i}(v)$.
The scalar coefficients of the $a_{g_i}$ may be varied independently
to obtain a vector space of dimension $n$ parametrizing the possible
twisted graded Hecke algebras realizable by the formula (\ref{eqn:udf})
and scalar modifications. It is shown in \cite[Example 2.16]{witherspoon}
that these are in fact {\em all} the twisted graded Hecke algebras
for this choice of $G$ and $\alpha$ in case $\ell\neq 2$.

In case $\ell =2$ and $n=3$, the elements $g_1,g_2,g_3$ are precisely
the nonidentity elements of $G$. The
formal deformation of $S(V)\#_{\alpha}G$ arising from the formula 
(\ref{eqn:udf}) is nearly the universal deformation, as is 
justified by considering the Hochschild cohomology of $S(V)\#_{\alpha}G$
(see \cite[Example 4.7]{caldararu-giaquinto-witherspoon} or the more
general Corollary \ref{cor:additive}).
That is, every Hochschild two-cocycle
$\mu_1$ with image in $S(V)\#_{\alpha}(G-\{1\})$ is an infinitesimal
of the formal deformation resulting from (\ref{eqn:udf})
with appropriate choices of the polynomials $s$ in (\ref{eqn:d2}).
(Classical deformations corresponding to the choice $g=1$
involve derivations that do not commute with the actions of the $H_i$,
and so we do not include these in the formula.) If $n>3$ or $\ell >2$,
there are nonidentity group elements other than $g_1,\ldots,g_n$, 
and the actions of the corresponding
Hopf algebras may no longer commute (but see the next example
below).
If $\alpha$ is not taken to be the cocycle we have chosen, the 
images of the $H_i$  in $\End_{\C}(S(V)\#_{\alpha}G)$ again
may no longer commute, and we do not know
whether there is a universal deformation formula
more complicated than (\ref{eqn:udf}) involving these operators.}
\end{example}

\begin{example}\label{ex:2}{\em
Let $G$ be a group acting on a vector space $V$ of dimension $n$,
$g$ a central element of $G$,
and assume (\ref{eqn:detg}) and (\ref{eqn:preserves}) hold with $q$
a primitive $\ell$th root of unity, $\ell >2$. Thus
$g$ corresponds to $\diag(q,q^{-1},1,\ldots,1)$
and $g^{-1}$ corresponds to $\diag(q^{-1},q,1,\ldots,1)$.
Assume further that $\alpha(g,g^{-1})=\alpha(g^{-1},g)$, as is true
in the last example for $g=g_i$. (In case $G$ is finite, this
assumption imposes no loss of generality, as any two-cocyle is
cohomologous to one satisfying this assumption 
\cite[Thm.\ 3.6.2]{karpilovsky85}).
Consider the images of $H_q$ and $H_{q^{-1}}$ in
$\End_{\C}(S(V)\#_{\alpha}G)$, 
where we let $D_1^{g^{-1}}$ involve $q^{-1}$-differentiation with respect
to $x_1$ and $D_2^{g^{-1}}$ involve $q$-differentiation with respect to 
$x_2$ in (\ref{eqn:d1}) and (\ref{eqn:d2}). 
Multiplying and dividing
the left side of the equation below by $q^{i-2}q^{i-1}$ yields
the right side:
$$
  \frac{(i)_q(i-1)_{q^{-1}}}{(i)_{q^{-1}}(i-1)_q} = q.
$$
Using this identity and the assumption $\alpha(g,g^{-1})=\alpha(g^{-1},g)$, 
we find that the following relations 
hold in $\End_{\C}(S(V)\#_{\alpha}G)$ among the images of the
generators of $H_q$ and $H_{q^{-1}}$:
\begin{eqnarray*}
  D_1^{g^{-1}} D_1^g = qD_1^gD_1^{g^{-1}} &,& 
         D_2^{g^{-1}}D_2^g = q^{-1}D_2^g D^{g^{-1}}_2,\\
  D_1^g D_2^{g^{-1}} = D_2^{g^{-1}}D_1^g &,&
         D_1^{g^{-1}}D_2^g = D_2^g D_1^{g^{-1}},\\
  \sigma^g\sigma^{g^{-1}} \!\! &=&\!\! \sigma^{g^{-1}}\sigma^g,\\
  \sigma^gD_1^{g^{-1}}=q^{-1}D_1^{g^{-1}}\sigma^g &,& \sigma^gD_2^{g^{-1}}
        =qD_2^{g^{-1}}\sigma^g,\\
  \sigma^{g^{-1}}D_1^g = qD_1^g\sigma^{g^{-1}} &,& 
       \sigma^{g^{-1}}D_2^g = q^{-1}D_2^g \sigma^{g^{-1}}.
\end{eqnarray*}
These relations are preserved by $\Delta$, $\varepsilon$ and $S$, 
and so the algebra generated by $H_q$ and $H_{q^{-1}}$, subject
to the above relations, is a Hopf algebra.
The proof of Theorem \ref{thm:main}
may be modified to show that $\exp_q(tD_1^g\ot D_2^g) \exp_{q^{-1}}
(tD_1^{g^{-1}}\ot D_2^{g^{-1}})$ is a universal deformation formula.
The key idea is to move factors corresponding to $g^{-1}$
past factors corresponding to $g$, so that the proof of Theorem 
\ref{thm:main} may be applied separately for each of $g,g^{-1}$. The relations
above imply that indeed the appropriate factors commute.
In case $G=\Z/3\Z$, this formula will nearly result in a 
universal deformation (again having infinitesimal with image in
$S(V)\#_{\alpha}(G-\{1\})$).}
\end{example}

 \bigskip

\section{Appendix: A computation of Hochschild cohomology}

The {\em Hochschild cohomology} of a $\C$-algebra $A$
is $\HH^{\DOT}(A):= \Ext^{\DOT}
_{A^e}(A,A)$, where $A^e=A\ot A^{op}$ acts on $A$ by left and
right multiplication. More generally, if $M$ is an $A$-bimodule
(equivalently, an $A^e$-module), we may define $\HH^{\DOT}(A,M):=
\Ext^{\DOT}_{A^e}(A,M)$, so that $\HH^{\DOT}(A)=\HH^{\DOT}(A,A)$.
These Ext groups may be expressed via the $A^e$-free resolution
of $A$: 
\begin{equation}\label{eqn:bar}
  \cdots\stackrel{\delta_3}{\longrightarrow}
   A^{\ot 4}\stackrel{\delta_2}{\longrightarrow}
  A^{\ot 3}\stackrel{\delta_1}{\longrightarrow}
  A^e \stackrel{m}{\longrightarrow} A\rightarrow 0,
\end{equation}
where $m$ is multiplication and
$$
   \delta_i(a_0\ot a_1\ot\cdots\ot a_{i+1}) =\sum_{j=0}^i (-1)^j
    a_0\ot\cdots\ot a_ja_{j+1}\ot\cdots\ot a_{i+1}.
$$
Applying $\Hom_{A^e}( - , M)$ and dropping the term 
$\Hom_{A^e}(A,M)$, we obtain
$$
0\rightarrow\Hom_{A^e}(A^e,M)\stackrel{\delta_1^*}{
  \longrightarrow} \Hom_{A^e}(A^{\ot 3},
  M) \stackrel{\delta_2^*}{
  \longrightarrow} \Hom_{A^e}(A^{\ot 4},M)\stackrel{\delta_3^*}
  {\longrightarrow} \cdots
$$
Then $\HH^i(A,M)=\Ker(\delta_{i+1}^*)/\Ima(\delta_i^*)$ and
$\HH^{\DOT}(A,M)=\oplus_{i\geq 0}\HH^i(A,M)$. Noting that
$\Hom_{A^e}(A^{\ot (i+2)},A)\cong \Hom_{\C}(A^{\ot i},A)$,
a straightforward calculation shows that $\HH^2(A)$ may be 
identified with the space of $\C$-linear functions $\mu_1:
A\ot A\rightarrow A$ satisfying the Hochschild two-cocycle
condition (\ref{eqn:cocycle}), modulo coboundaries.
See \cite{weibel94} for more details on Hochschild cohomology.

Let $G$ be a finite subgroup of ${\rm {GL}}(V)$. We will
compute $\HH^{\DOT}(S(V)\#_{\alpha}G)$, using techniques similar to those
in \cite{AFLS,suarez-alvarez}, 
where the crossed product was taken with a Weyl algebra
rather than a polynomial algebra. 
We will use a result of \c Stefan on Hopf Galois extensions 
\cite[Cor.\ 3.4]{stefan95}. It implies that there is an action
of $G$ on $\HH^{\DOT}(S(V),S(V)\#_{\alpha}G)$ for which
\begin{equation}\label{eqn:stefan}
  \HH^{\DOT}(S(V)\#_{\alpha}G)\cong \HH^{\DOT}(S(V),
  S(V)\#_{\alpha}G)^G,
\end{equation}
where the superscript $G$ denotes the subspace of 
$G$-invariant elements. (A more explicit proof of this result,
useful in this context,
is given in \cite[\S5]{caldararu-giaquinto-witherspoon}).
A Koszul complex may then be used to compute $\HH^{\DOT}
(S(V),S(V)\#_{\alpha}G)$. This is done in a more general geometric setting
by Ginzburg and Kaledin 
\cite{ginzburg-kaledin04} in the case $G$ is symplectic
and $\alpha$ is trivial, although they note that their techniques
apply to any finite group $G$. 
An elegant algebraic computation is given by Farinati \cite{farinati}
for an arbitrary finite group $G$, and trivial $\alpha$.
The additive structure of
$\HH^{\DOT}(S(V),S(V)\#_{\alpha}G)$, before taking $G$-invariants, 
is independent of $\alpha$
since the $S(V)$-bimodule structure of $S(V)\#_{\alpha}G$
does not involve $\alpha$. 
Thus the techniques of either \cite{farinati} or \cite{ginzburg-kaledin04}
apply here. For completeness, we give an explicit algebraic computation
whose details are needed in Section 4.

Note that $S(V)\#_{\alpha}G=\oplus_{g\in G}S(V)\overline{g}$, 
where $S(V)\overline{g} =\{s\overline{g}\mid s\in S(V)\}$, as an 
$S(V)$-bimodule. 
Thus there is an additive 
decomposition of Hochschild cohomology,
\begin{equation}\label{eqn:directsum}
  \HHD(S(V),S(V)\#_{\alpha}G)\cong \bigoplus_{g\in G}\HHD(S(V),
    S(V)\overline{g}).
\end{equation}
We will determine each summand $\HHD(S(V),S(V)\overline{g})$, noting again
that $\alpha$ plays
no role here as we need only the $S(V)$-module structure of each
$S(V)\overline{g}$.
If $g=1$, we have $\HHD(S(V),S(V)\overline{1})=\HHD(S(V))$, 
and the Hochschild-Kostant-Rosenberg
Theorem states that
$$
  \HHD(S(V))\cong\Wedge^{\DOT}_{S(V)}(S(V)^n)\cong
  \Wedge^{\DOT}(V^*)\ot S(V),
$$
where $n=\dim V$.
(See \cite[Exer.\ 9.1.3 and Thm.\ 9.4.7]{weibel94}.) 
Letting $x_1,\ldots,x_n$ be a basis of $V$,
this may be computed directly from the $S(V)^e$-projective Koszul resolution 
$K(\{x_i\ot 1-1\ot x_i\}_{i=1}^n)\cong \Wedge^{\DOT}(V)\ot S(V)^e$ 
of $S(V)$. (See \cite[\S4.5]{weibel94} for details on Koszul complexes.)
The differential $d_m:\Wedge ^m(V)\ot S(V)^e$ is given by
$$
  d_m(e_{i_1}\wedge\cdots\wedge e_{i_m} \ot 1\ot 1) =
  \sum_{k=1}^m (-1)^{k+1} e_{i_1}\wedge \cdots\wedge \widehat{e_{i_k}}\wedge
  \cdots\wedge e_{i_m}\ot (x_{i_k}\ot 1 - 1\ot x_{i_k}),
$$
where we use the standard notation $e_{i_j}$ for the element $x_{i_j}$ in
$\Wedge^1(V)$.
After application of $\Hom_{S(V)^e}
(-,S(V))$, all chain maps become 0 as $S(V)$ is commutative. 
We will identify $\HHD(S(V))$ with $\Wedge^{\DOT}(V^*)\ot S(V)$, 
as the group action is clear
in that notation: It is diagonal on the factors, 
with the standard actions on $S(V)$ and
on $\Wedge^{\DOT}(V^*)$. In case of an element $g$ not necessarily
equal to 1, we have the following.

\begin{prop}\label{thm:additive}
For each $g\in G$, 
$$\HHD(S(V),S(V)\overline{g})\cong
\Wedge^{\DOT-\codim V^g}((V^g)^*)\ot S(V^g).$$
\end{prop}

In particular, the lowest degree $j$ for which $\HH^{j}(S(V),S(V)\overline{g})
\neq 0$ is $j=\codim V^g$.

\begin{proof} Fix $g\in G$. As the order of $g$ is finite, we may assume
without loss of generality that the action of $g$ is diagonal with respect
to the basis $x_1,\ldots,x_n$ of $V$. 
Thus there are scalars 
$\lambda_i$ with $g\cdot x_i=\lambda_i x_i \ (i=1,\ldots,n)$. We will
further assume, for notational convenience, that the basis is ordered
so that $\lambda_i=1$ for $i=1,\ldots,r$ and $\lambda_i\neq 1$ for 
$i=r+1,\ldots,n$. (Set $r=0$ if $V^g=0$ and $r=n$ if $V^g=V$.) 
We may also assume that $\Span_{\C}\{x_{r+1},\ldots,x_n\}=(V^g)^{\perp}$
where the orthogonal complement is taken with respect to some nondegenerate
$G$-invariant Hermitian form on $V$.

Consider the complex $\Hom_{S(V)^e}(K(\{x_i\otimes 1 - 1\ot x_i\}_{i=1}^n), 
S(V)\overline{g})$,
which we may identify with $\Hom_{S(V)^e}(\Wedge^{\DOT}(V)\ot S(V)^e,S(V)
\overline{g})\cong\Wedge^{\DOT}(V^*)\ot S(V)\overline{g}$. 
Additively, this is the same as $\Wedge^{\DOT}(V^*)\ot S(V)$, 
but the factor $\overline{g}$ affects the differentials, which
we will determine next.
They are {\em not} necessarily all zero (in contrast to the case $g=1$).
If $s\in S(V)$, we have
$$
  (x_i\ot 1 - 1\ot x_i)\cdot s\overline{g} 
  = x_is\overline{g}-s\overline{g}x_i = (x_i-g\cdot x_i)s\overline{g}.
$$
If $i=1,\ldots,r$, this element is 0. If $i=r+1,\ldots,n$,
the factor $x_i-g\cdot x_i=(1-\lambda_i)x_i$ is a nonzero scalar multiple 
of $x_i$.
Thus $\Hom_{S(V)^e}(K(\{x_i\ot 1-1\ot x_i\}_{i=1}^n),S(V)
\overline{g})$ is equivalent
to the dual Koszul complex $\overline{K}(0,\ldots,0,x_{r+1},\ldots, x_n)$
for $S(V)$, where the bar denotes the reverse order. This is the tensor
product (over $S(V)$) of two complexes for $S(V)$: $\overline{K}(0,\ldots,0)$
and $\overline{K}(x_{r+1},\ldots,x_n)$. The second complex is exact 
other than in degree $n-r$ (as the corresponding Koszul complex is
exact other than in degree 0), where it has cohomology
$S(V)/(x_{r+1},\ldots,x_n)S(V)\cong S(V^g)$ (see \cite[Cor.\ 4.5.4]{weibel94}).
We will identify this with $\det(((V^g)^{\perp})^*)\ot S(V^g)$, where 
$\det(((V^g)^{\perp})^*)$ is the one-dimensional space 
$\Wedge^{\codim V^g} (((V^g)^{\perp})^*)$, to account for the degree shift
and the action of $G$.
The spectral sequence of
the double complex $\overline{K}(0,\ldots,0)\ot_{S(V)}\overline{K}(x_{r+1},
\ldots,x_n)$ thus collapses at $E_2$ with $E_2^{pq}=0$ for $q\neq n-r$,
and $E_2^{p,\ n-r} ={\rm {H}}^p(\overline{K}(0,\ldots,0))\ot_{S(V)} S(V^g)$
by freeness of the terms of the chain complex over $S(V)$.
This follows from \cite[Thm.\ 3.4.2]{benson91b}, which also implies 
that the cohomology is precisely $E_2^{p, \ n-r}$.
Now ${\rm{H}}^p(\overline{K}(0,\ldots,0))\cong\Wedge^p_{S(V)}(S(V)^{\dim V^g})
\cong \Wedge^p((V^g)^*)\ot S(V)$,
and as this is tensored with the cohomology of $\overline{K}(x_{r+1},\ldots,
x_n)$, namely $S(V)/(x_{r+1},\ldots,x_n)S(V)\cong S(V^g)$ in degree $n-r$, 
we obtain the stated result.
\end{proof}

We will identify the cohomology 
$\Wedge^{\DOT-\codim V^g}((V^g)^*)\ot S(V^g)$ of the theorem 
with $\Wedge^{\DOT -\codim V^g}((V^g)^*)\ot \det (((V^g)^{\pe})^*)
\ot S(V^g)\overline{g}$. 
The action of $G$ is nontrivial on the one-dimensional 
factor $\det(((V^g)^{\pe})^*)$, as may be seen by considering the
action of $G$ on the corresponding cochain complex, and the action
on $\overline{g}$ is by conjugation by $\overline{h}$ ($h\in G$).

The following corollary is immediate from (\ref{eqn:stefan}), 
(\ref{eqn:directsum}) and Proposition \ref{thm:additive}, after
making the above identifications.

\begin{cor}\label{cor:additive}
There is an additive decomposition of Hochschild cohomology,
$$
  \HHD(S(V)\#_{\alpha}G)\cong \left(\bigoplus_{g\in G} 
  \Wedge^{\DOT-\codim V^g}((V^g)^*)\ot\det(((V^g)^{\perp})^*)\ot S(V^g)
   \overline{g}\right)^G.
$$
\end{cor}

Compare the above corollary with \cite[Thm.\ 3.6]{farinati}, or
with the formula just above (6.4) in \cite{ginzburg-kaledin04}.

As a consequence of Corollary \ref{cor:additive}, we obtain a necessary
condition for there to exist a Hochschild
two-cocycle $\mu_1$, with image in the $g$-component 
$S(V)\overline{g}$, that is not a coboundary.
Due to the degree shift $2\! -\!\codim V^g$, such an element $g$ must satisfy
$\codim V^g\in \{0,1,2\}$. 
We claim that the determinant of $g$ on $V$ must be 1.
If $\det(g)\neq 1$, the action of $g$ itself
on the one-dimensional space $\det(((V^g)^{\perp})^*)$ is nontrivial, 
whereas its actions on $\Wedge^{2-\codim V^g}((V^g)^*)$ 
and on $S(V^g)\overline{g}$ are trivial.
Consequently there can be no such $G$-invariant elements.
Therefore $\det(g)=1$,
which also now implies $\codim V^g\in \{0,2\}$.
(See also \cite[Ex.\ 3.10]{farinati}.)


\begin{thebibliography}{99}

\bibitem{AFLS} J.\ Alev, M.\ A.\ Farinati, T.\ Lambre, and A.\ L.\
Solotar, ``Homologie des invariants d'une alg\`{e}bre de Weyl sous
l'action d'un groupe fini,'' J.\ Algebra 232 (2000), 564--577.

\bibitem{benkart-witherspoon} G.\ Benkart and S.\ Witherspoon,
``Quantum group actions, twisting elements, and deformations of algebras,''
{\tt math.QA/0503310}.

\bibitem{benson91b} D.\ J.\ Benson, {\it Representations and Cohomology
II: Cohomology of Groups and Modules}, Cambridge Univ.\ Press,
Cambridge, 1991.

\bibitem{caldararu-giaquinto-witherspoon} A.\ C\u ald\u araru, A.\ Giaquinto,
and S.\ Witherspoon, ``Algebraic deformations arising from orbifolds with
discrete torsion,'' J.\ Pure Appl.\ Algebra 187 (2004), no.\ 1--3, 51--70.

\bibitem{cibils90} C.\ Cibils, ``Rigidity of truncated quiver
algebras,'' Adv.\ Math.\ 79 (1990), 18--42.

\bibitem{cibils91} C.\ Cibils, ``Rigid monomial algebras,'' Math.\ Ann.\
289 (1991), 95--109.

\bibitem{drinfeld86} V.\ G.\ Drinfel'd, ``Degenerate affine Hecke
algebras and Yangians,'' Funct.\ Anal.\ Appl.\ 20 (1986), 58--60.

\bibitem{etingof-ginzburg} P.\ Etingof and V.\ Ginzburg, ``Symplectic 
reflection algebras, Calogero-Moser space, and deformed Harish-Chandra
homomorphisms,'' Invent.\ Math.\ 147 (2002), no.\ 2, 243--348.

\bibitem{farinati} M.\ Farinati, ``Hochschild duality, localization,
and smash products,'' J.\ Algebra 284 (2005), no.\ 1, 415--434.

\bibitem{gerstenhaber-schack92} M.\ Gerstenhaber and S.\ D.\ Schack,
``Algebras, bialgebras, quantum groups, and algebraic deformations,''
in {\it Deformation Theory and Quantum Groups with Applications
to Mathematical Physics}, Contemp.\ Math.\ 134 (1992), 51--92.

\bibitem{giaquinto-zhang98} A.\ Giaquinto and J.\ Zhang, ``Bialgebra
actions, twists, and universal deformation formulas,''  J.\ Pure
Appl.\ Algebra 128 (1998), 133--151.

\bibitem{GZ} A.\ Giaquinto and J.\ Zhang, private communication.

\bibitem{ginzburg-kaledin04} V.\ Ginzburg and D.\ Kaledin, ``Poisson deformations
of symplectic quotient singularities,'' Adv.\ Math.\ 186 (2004), no.\ 1, 1--57.

\bibitem{karpilovsky85} G.\ Karpilovsky, {\em Projective representations
of finite groups}, Marcel-Dekker, 1985.

\bibitem{kassel95} C.\ Kassel, {\it Quantum Groups}, Springer-Verlag, 1995.

\bibitem{kharchenko91} V.\ K.\ Kharchenko, {\it Automorphisms and
Derivations of Associative Rings}, Kluwer, 1991.

\bibitem{kharchenko01} V.\ K.\ Kharchenko, ``Skew primitive elements
in Hopf algebras and related identities,'' J.\ Algebra 238 (2001), 534--559.

\bibitem{lusztig89} G.\ Lusztig, ``Affine Hecke algebras and their
graded version,'' J.\ Amer.\ Math.\ Soc.\ 2 (1989), 599--635.

\bibitem{montgomery93} S.\ Montgomery, {\it Hopf Algebras and Their
Actions on Rings}, CBMS Regional Conference Series in Math., Number 82,
Amer.\ Math.\ Soc., 1993.

\bibitem{passman89} D.\ S.\ Passman, {\it Infinite Crossed Products},
Academic Press, New York, 1989.

\bibitem{ram-shepler} A.\ Ram and A.\ V.\ Shepler, ``Classification of
graded Hecke algebras for complex reflection groups,'' Comment.\ Math.\ Helv.\
78 (2003), no.\ 2, 308--334.

\bibitem{stefan95} D.\ \c{S}tefan, ``Hochschild cohomology on Hopf Galois
extensions,'' J.\ Pure Appl.\ Algebra 103 (1995), 221--233.

\bibitem{suarez-alvarez} M.\ Suarez-Alvarez, ``Algebra structure on the
Hochschild cohomology of the ring of invariants 
of a Weyl algebra under a finite group,'' J.\ Algebra
248 (2002), 291--306.

\bibitem{weibel94} C.\ A.\ Weibel, {\it An Introduction to Homological
Algebra}, Cambridge Studies in Adv.\ Math.\ 38, Cambridge Univ.\ Press,
Cambridge, 1994.

\bibitem{witherspoon} S.\ Witherspoon, ``Twisted graded Hecke algebras,''
{\tt math.RT/0506152}.

\end{thebibliography}
\end{document}